\documentclass[11pt]{article}
   \usepackage{amsfonts}
\usepackage{latexsym,tikz}
\usepackage{amsmath}
 
\title{Multiple Dedekind Zeta Functions}
\author{Ivan Emilov Horozov}
\date{February 2, 2008}

\newcommand \nc {\newcommand}
\nc \proof {\noindent {\em{Proof.\/ }}} \nc \qed {$\Box$\hfill}
\newtheorem{theorem}{Theorem}[section]
\newtheorem{lemma}[theorem]{Lemma}
\newtheorem{proposition}[theorem]{Proposition}
\newtheorem{corollary}[theorem]{Corollary}
\newtheorem{definition}[theorem]{Definition}
\newtheorem{example}[theorem]{Example}
\newtheorem{remark}[theorem]{Remark}
\newtheorem{conjecture}[theorem]{Conjecture}
\newtheorem{question}[theorem]{Question}
\nc \bth[1] {\begin{theorem}\label{t#1} } \nc \ble[1]
{\begin{lemma}\label{l#1} } \nc \bpr[1]
{\begin{proposition}\label{p#1} } \nc \bco[1]
{\begin{corollary}\label{c#1} } \nc \bde[1]
{\begin{definition}\label{d#1}\rm } \nc \bex[1]
{\begin{example}\label{e#1}\rm } \nc \bre[1]
{\begin{remark}\label{r#1}\rm } \nc \bcon[1]
{\begin{conjecture}\label{con#1}\rm } \nc \bque[1]
{\begin{question}\label{que#1}\rm }
\nc {\eth} { \end{theorem} } \nc {\ele} { \end{lemma} } \nc
{\epr}{ \end{proposition} } \nc {\eco} { \end{corollary} } \nc
{\ede} {\end{definition} } \nc {\eex} { \end{example} } \nc {\ere}
{\end{remark} } \nc {\econ} { \end{conjecture} } \nc {\eque}
{\end{question} }

\def \N {{\mathbb N}}

\def \Q {{\mathbb Q}}
\def \R {{\mathbb R}}

\begin{document}

\title{{\LARGE\bf{Double Shuffle Relations\\ for Multiple Dedekind Zeta Values}}}

\author{
Ivan ~Horozov
\thanks{E-mail: horozov@math.wustl.edu}
\\ \hfill\\ \normalsize \textit{Department of Mathematics,}\\
\normalsize \textit{Washington University in St. Louis,}\\
\normalsize \textit{One Brookings Drive, Campus box 1146}\\
\normalsize \textit {Saint Louis, MO 63130, USA }  \\
}
\date{}
\maketitle

\begin{abstract}
This paper contains examples of shuffle relations among multiple Dedekind zeta values.
Dedekind zeta values were defined by the author in \cite{MDZF}.
Here we concentrate on the cases of real or imaginary quadratic fields with up to double iteration.
We give examples of integral shuffle relation in terms of iterated integrals over membranes and of infinite sum shuffle relation, sometimes called stuffle relation. Using both types of shuffles for the product $\zeta_{K,C}(2)\zeta_{K,C}(2)$, we find relations among multiple Dedekind zeta values 
for real and for imaginary quadratic fields.
\end{abstract}

\tableofcontents
\setcounter{section}{-1}

\section{Introduction}
\label{Intro}

Multiple zeta values were defined by Euler \cite{Eu}:
\begin{equation}
\label{eq MZV}
\zeta(m_1,\dots,m_d)=\sum_{0<k_1<\dots<k_d}\frac{1}{{k_1}^{m_1}\dots {k_d}^{m_d}}.
\end{equation}
They can be expressed as iterated path integrals. This integral representation leads to an integrals shuffle relation among the multiple zeta values. It gives a formula expressing a product of multiple zeta values as a sum of such. Using the infinite sum representation \eqref{eq MZV}, again, we can express a product of two multiple zeta values as a sum of such. However, the infinite sum representation gives different formulas for the product compared to the shuffle coming from the integral representation. Expressing a product of multiple zeta values as a sum of such in two different ways leads to linear relations among multiple zeta values (see for example \cite{G} \cite{G2}, \cite{GKZ}).

Multiple Dedekind zeta values are defined in \cite{MDZF} using a higher dimensional analogue of iterated path integrals, which we call iterated integrals over membranes. The idea of iterated integrals over membrane is written first by the author in \cite{ModSym} for the purpose of generalizing Manin's non-commutative modular symbol \cite{Man}  to Hilbert modular surfaces. That project is completed in \cite{Hilbert}.

Similarly to multiple zeta values, multiple Dedekind zeta values also have two shuffle relations. Using the representation of multiple Dedekind zeta values in terms of iterated integrals over membranes, one obtain a shuffle relation, expressing a product of multiple Dedekind zeta values as a sum of such. Examples of this type are presented in this paper. Also in this paper, we give examples of infinite sum shuffle relation for real and for imaginary quadratic fields. The formulas for real quadratic fields can easily be generalized to any totally real number field. Using both types of shuffles, we obtain relations among multiple Dedekind zeta values.

For imaginary quadratic fields, we use (in the last Section) an idea of Gangl, Kaneko, Zagier \cite{GKZ}, for separating the ring of integers in an imaginary quadratic field into a positive region, a negative region and zero. 

In Section 1, we give the definitions of multiple Dedekind zeta values needed for the two types of shuffle.

In Section 2, we recall the infinite sum shuffle for multiple zeta values. We construct analogous infinite sum shuffle for multiple Dedekind zeta values. In particular, we obtain an infinite sum shuffle for a product of two Dedekind zeta values at $s=2$. This shuffle relation expresses a product of two Dedekind zeta values at $s=2$ as a sum of multiple Dedkind zeta values.

In Section 3, we recall classical polylogarithms. We express (multiple) polylogarithms as iterated integrals.
In particular, that gives multiple zeta values in the form of iterated integrals. We give examples of integral shuffle of multiple zeta values. Also, we express the polylogarithms as integrals involving exponent, since this idea is generalized. We call this generalization multiple Dedekind polylogarithms in the form of iterated integrals over membranes. In particular, that expresses multiple Dedekind zeta values as iterated integrals over membranes. There is an example of integral shuffle relation for Dedekind polylogarithms both in terms of iterated integrals. We also represent such integrals in terms of diagrams. At the end of the section, consider an integral shuffle relation of two Dedekind zeta values at $s=2$. This shuffle relation expresses the product as a sum of multiple Dedekind zeta values. For convenience of the reader, we have included all diagrams that appear in such shuffle and the corresponding terms that the diagrams represent. 

In Section 4, we use that the integral shuffle and the infinite sum shuffle of two Dedekind zeta values at $s=2$ give different formulas in terms of sums of multiple Dedekind zeta values. Expressing a product  as sums in two different ways gives us a linear relation among multiple Dedekind zeta values.
For real and for imaginary quadratic fields that produces different formulas due to the fact that the infinite sum shuffles for the real and for imaginary quadratic fields are different, even though the integral shuffle for both types of fields are, in principle, the same.

\section{Examples of multiple Dedekind zeta values}
In this section we give many examples of multiple Dedekind zeta values from the point of view of shuffle relations.

First, let us recall double zeta value
\[\zeta(a,b)=\sum_{0<m<n}\frac{1}{m^an^b}.\]
The double zeta values can also be written as
\[\zeta(a,b)=\sum_{m,n=1}^\infty\frac{1}{m^a(m+n)^b}.\]

Now let us examine a real quadratic field $K$. Let ${\cal{O}}_K$ be the ring of integers in the number field $K$. For example, the elements of the field $K=\Q(\sqrt{2})$ are all numbers of the form 
$p+q\sqrt{2}$, where $p$ and $q$ are rational numbers. 
Then, the algebraic integers in the number field $K=\Q(\sqrt{2})$ are of the form $m+n\sqrt{2}$, where $m$ and $n$ are integers.
If $\alpha_1=\alpha=m+n\sqrt{d}$ then its Galois conjugate is $\alpha_2=m-n\sqrt{d}$.
We call an element $\alpha\in K$ totally positive if both $\alpha_1=\alpha$ and its Galois conjugate 
$\alpha_2$ are positive real numbers. For example, $5+\sqrt{2}$ is totally positive, because both
$5+\sqrt{2}>0$ and its Galois conjugate $5-\sqrt{2}>0$  are positive.
However, $1+\sqrt{2}$ is not totally positive, because
its Galois conjugate is negative, $1-\sqrt{2}<0$.

Let $\mu$ and $\nu$ be totally positive integers in the real quadratic number field $K$.
We define a cone as a subset of the algebraic integers in $K$ 
generated by $\mu$ and $\nu$ over $\N$, namely,
\[C=\N\{\mu,\nu\}=\{\alpha\in {\cal{O}}_K\,|\,\alpha=m\mu+n\nu\}.\]

In the paper \cite{MDZF}, we have defined multiple Dedekind zeta values of the form:
 \[\zeta_{K,C}(a)=\sum_{\alpha\in C}\frac{1}{N(\alpha)^a}\]
and
 \[\zeta_{K,C}(a,b)=\sum_{\alpha,\beta\in C}\frac{1}{N(\alpha)^aN(\alpha+\beta)^b}.\]
as terms of iterated integrals over membranes, where $N(\alpha)=\alpha_1\alpha_2$. Using the same technique, we are lead to the definition of more complex infinite sums, which are needed for the purpose of shuffle relations.

\begin{definition} For the purpose of establishing integral shuffle relation, we define the following multiple Dedekind zeta values:
\[\zeta^1_{K,C}(a_1,b_1;c_2,d_2)=
\sum_{\alpha,\beta\in C}
\frac{1}{\alpha_1^{a_1}(\alpha_1+\beta_1)^{b_1}\alpha_2^{c_2}(\alpha_2+\beta_2)^{d_2}}\]
and
\[\zeta^\rho_{K,C}(a_1,b_1;c_2,d_2)=
\sum_{\alpha,\beta\in C}
\frac{1}{\beta_1^{a_1}(\alpha_1+\beta_1)^{b_1}\alpha_2^{c_2}(\alpha_2+\beta_2)^{d_2}}.\]
In particular
\[\zeta^1_{K,C}(a,b;a,b)=\zeta_{K,C}(a,b).\]
\end{definition}

In order to have an integral shuffle of a product $\zeta_{K,C}(2)\zeta_{K,C}(2)$, we need multiple Dedekind zeta values of the types $\zeta^1$ and $\zeta^\rho$. Such relation will be established in Section 3.

In order to have infinite sum shuffle we need the following 

\begin{definition}
Let $\mu(\alpha)=m$ and  $\nu(\alpha)=n$ if $\alpha=m\mu+n\nu$ for integers values of $m$ and $n$.
The following multiple Dedekind zeta values are needed for the infinite sum shuffle of Dedekind zeta values:
\[_1\zeta_{K,C}(a_1,b_2;c_2,d_2)=\sum_{\alpha,\beta\in C\,:\,\mu(\alpha)<\mu(\beta),\,\nu(\alpha)<\nu(\beta)}
\frac{1}{\alpha_1^{a_1}\beta_1^{b_1}\alpha_2^{c_2}\beta_2^{d_2}},\]
and
\[_\rho\zeta_{K,C}(a_1,b_2;c_2,d_2)=\sum_{\alpha,\beta\in C\,:\,\mu(\alpha)>\mu(\beta),\,\nu(\alpha)<\nu(\beta)}
\frac{1}{\alpha_1^{a_1}\beta_1^{b_1}\alpha_2^{c_2}\beta_2^{d_2}}.\]
We also need a couple of boundary Dedekind zeta values in order to have an infinite sum shuffles. Let us define
\[_{0,1}\zeta_{K,C}(a_1,b_2;c_2,d_2)=
\sum_{\alpha,\beta\in C\,:\,\mu(\alpha)=\mu(\beta),\,\nu(\alpha)<\nu(\beta)}
\frac{1}{\alpha_1^{a_1}\beta_1^{b_1}\alpha_2^{c_2}\beta_2^{d_2}}\]
and
\[_{1,0}\zeta_{K,C}(a_1,b_2;c_2,d_2)=
\sum_{\alpha,\beta\in C\,:\,\mu(\alpha)<\mu(\beta),\,\nu(\alpha)=\nu(\beta)}
\frac{1}{\alpha_1^{a_1}\beta_1^{b_1}\alpha_2^{c_2}\beta_2^{d_2}}.\]
\end{definition}

\begin{lemma}
We have the following relation among multiple Dedekind zeta values
\[_1\zeta(a_1,b_2;c_2,d_2)=\zeta(a_1,b_2;c_2,d_2)^1.\]
\end{lemma}
\proof Let $\alpha$ and $\beta$ be in $C$ so  that
\[\mu(\alpha)<\mu(\beta),\,\nu(\alpha)<\nu(\beta)\] Let $\gamma=\beta-\alpha$. 
Then $\mu(\gamma)=\mu(\beta)-\mu(\alpha)>0$ and 
$\mu(\gamma)=\nu(\beta)-\nu(\alpha)>0$. Therefore $\gamma\in C$. Conversely,
if $\gamma\in C$ and $\alpha\in C$ then for $\beta=\alpha+\gamma$ we have $\beta\in C$. Therefore, $\mu(\beta>0$ and $\nu(\beta)>0$.
Then for $\gamma=\beta-\alpha$, we have
\begin{align*}
_1\zeta_{K,C}(a_1,b_2;c_2,d_2)
=&\sum_{\alpha,\beta\in C\,:\,\mu(\alpha)<\mu(\beta),\,\nu(\alpha)<\nu(\beta)}
\frac{1}{\alpha_1^{a_1}\beta_1^{b_1}\alpha_2^{c_2}\beta_2^{d_2}}=\\
=&\sum_{\alpha,\gamma\in C}
\frac{1}{\alpha_1^{a_1}(\alpha_1+\gamma_1)^{b_1}\alpha_2^{c_2}(\alpha_2+\gamma_2)^{d_2}}=\\
=&\zeta_{K,C}(a_1,b_2;c_2,d_2)^1.
\end{align*}
\qed

The infinite sum shuffle relation expresses a product of two $\zeta_{K,C}(a)$ as a sum of 
the types $\,\,_1\zeta_{K,C}$, $\,\,_{0,1}\zeta_{K,C}$, $\,\,_{1,0}\zeta_{K,C}$ and $\,\,_\rho\zeta_{K,C}$.
The integral shuffle relation expresses a product of two $\zeta_{K,C}(a)$ as a sum of the types 
$\zeta^1_{K,C}$ and $\zeta^\rho_{K,C}$.

For imaginary quadratic fields $K$ the setting is simpler. 
Again, let ${\cal{O}}_K$ be the ring of integers in the number field $K$.
If $\alpha \in {\cal{O}}_K$ we define $\alpha_1=\alpha$ and $\alpha_2$ is the complex conjugate of 
$\alpha_1$. We will consider the norm of $\alpha$, $N(\alpha)$, which is defined as 
$N(\alpha)=\alpha_1\alpha_2$.
Following Gangl, Kaneko, Zagier, \cite{GKZ}, we define a cone $C_+$ by
\[C_+=\N\cup\{\alpha\in {\cal{O}}_K\,|\, Im(\alpha_1)>0\}.\]
Let also
\[C_-=\{-\alpha\,|\,\alpha\in C_+\}.\]
Then, we have a decomposition of the ring of integers into a finite disjoint union of sets
\[{\cal{O}}_K=C_+\cup\{0\}\cup C_-.\]

For imaginary quadratic fields, with the above definition of the cone $C_+$,  we  consider the following multiple Dedekind zeta values

\[\zeta^1_{K,C}(a_1,b_1;c_2,d_2)=
\sum_{\alpha,\beta\in C_+}
\frac{1}{\alpha_1^{a_1}(\alpha_1+\beta_1)^{b_1}\alpha_2^{c_2}(\alpha_2+\beta_2)^{d_2}}\]
and 
\[\zeta^\rho_{K,C}(a_1,b_1;c_2,d_2)=
\sum_{\alpha,\beta\in C_+}
\frac{1}{\beta_1^{a_1}(\alpha_1+\beta_1)^{b_1}\alpha_2^{c_2}(\alpha_2+\beta_2)^{d_2}}.\]

The infinite sum shuffle between two values of the type $\zeta^1_{K,C}(a,b)$ will be represented as double Dedekind zetas of the type $\zeta^1$. And the integral shuffle relation will represent a product of two  values of the type $\zeta^1_{K,C}(a,b)$ as a sum of double Dedekind zeta values of the types 
$\zeta^1$ and $\zeta^\rho$.



\section{Infinite sum shuffle relation}
We are going to recall the classical multiple zeta values and their shuffle relations.
The Riemann zeta values are defined as
\[\zeta(m)=\sum_{k=1}^\infty\frac{1}{k^m}.\]
Multiple zeta values are defined as
\[\zeta(m_1,\dots,m_d)=\sum_{0<k_1<\dots<k_d}\frac{1}{{k_1}^{m_1}\dots {k_d}^{m_d}}.\]
For multiple zeta values we have the following  sum shuffle relation, sometimes called stuffle relation:
\begin{align}
\zeta(m_1)\zeta(m_2)
&=\sum_{k_1=0}^\infty\frac{1}{{k_1}^{m_1}}\sum_{k_2=1}^\infty\frac{1}{{k_2}^{m_2}}=\\
&=\left(\sum_{0<k_1<k_2}+\sum_{0<k_1=k_2}+\sum_{0<k_2<k_1}\right)\frac{1}{{k_1}^{m_1}{k_2}^{m_2}}=\\
&=\zeta(m_1,m_2)+\zeta(m_1+m_2)+\zeta(m_2,m_1).
\end{align}
Similarly, a product of a zeta value and a double zeta value, can be expressed as a sum of multiple zeta values:
\begin{align}
\zeta(m_1)\zeta(m_2,m_3)
=&\sum_{k_1=1}^\infty\frac{1}{{k_1}^{m_1}}\sum_{k_2<m_3}^\infty\frac{1}{{k_2}^{m_2}{k_3}^{m_3}}=\\
=&\left(\sum_{0<k_1<k_2<m_3}+\sum_{0<k_1=k_2<m_3}+\sum_{0<k_2<k_1<k_3}+\right.\\
&+\left.\sum_{0<k_2<k_1=k_3}+\sum_{0<k_2<k_3<k_1}\right)\frac{1}{{k_1}^{m_1}{k_2}^{m_2}{k_3}^{m_3}}=\\
=&\zeta(m_1,m_2,m_3)+\zeta(m_1+m_2,m_3)+\zeta(m_2,m_1,m_3)+\\
&=\zeta(m_2,m_1+m_3)+\zeta(m_2,m_3,m_1).
\end{align}

Let $K$ be an imaginary quadratic field. Let ${\cal{O}}_K$ be the ring of integers in the field $K$. We are going to use the notation 
\[C_+=\N\cup\{\alpha\in {\cal{O}}_K\,|\, Im(\alpha_1)>0\}.\]
and
\[C_-=\{-\alpha\,|\,\alpha\in C_+\}.\]
Recall, the following decomposition
\[{\cal{O}}_K=C_+\cup\{0\}\cup C_-.\]
We write $\alpha<\beta$ for $\alpha$ and $\beta$ in $C_+$ if $\beta-\alpha$ is in $C_+$.
Let
\[\zeta^1_{K,C}(a_1;c_2)=
\sum_{\alpha,\beta\in C_+}
\frac{1}{\alpha_1^{a_1}\alpha_2^{c_2}}.\]
Then we define the following infinite sum shuffle relation
\begin{align*}
\zeta^1_{K,C}(a_1;c_2)\zeta^1_{K,C}(b_1;d_2)
=&\sum_{\alpha,\beta\in C_+}
\frac{1}{\alpha_1^{a_1}\alpha_2^{c_2}\beta_1^{b_1}\beta_2^{d_2}}=\\
=&\left(\sum_{\alpha<\beta}+\sum_{\alpha=\beta}+\sum_{\beta<\alpha}\right)
\frac{1}{\alpha_1^{a_1}\beta_1^{b_1}\alpha_2^{c_2}\beta_2^{d_2}}=\\
=&\zeta^1_{K,C}(a_1,b_1;c_2,d_2)
+\zeta^1_{K,C}(a_1+b_1;c_2+d_2)+\zeta^1_{K,C}(b_1,a_1;d_2,c_2)
\end{align*}
In particular,
\begin{equation}
\label{eq stuffle 2,2 C}
\zeta^1_{K,C}(2;2)\zeta^1_{K,C}(2;2)=\zeta^1_{K,C}(4;4)+2\zeta^1_{K,C}(2,2;2,2)
\end{equation}

Similarly, we define the following infinite sum shuffle
\begin{align*}
\zeta^1_{K,C}(a_1,b_1;d_2,e_2)\zeta^1_{K,C}(c_1;f_2)
=&\sum_{\alpha,\beta,\gamma\in C_+;\, \alpha<\beta}
\frac{1}{\alpha_1^{a_1}\beta_1^{b_1}\gamma_1^{c_2}\alpha_2^{d_2}\beta_2^{e_2}\gamma_2^{f_2}}=\\
=&\left(\sum_{\alpha<\beta<\gamma}+\sum_{\alpha<\beta=\gamma}+\sum_{\alpha<\gamma<\beta}+\right.\\
&\left. +\sum_{\gamma=\alpha<\beta}+\sum_{\gamma<\alpha<\beta}
\right)\frac{1}{\alpha_1^{a_1}\beta_1^{b_1}\gamma_1^{c_2}\alpha_2^{d_2}\beta_2^{e_2}\gamma_2^{f_2}}
=\\
=&\zeta^1(a_1,b_1,c_1;d_2,e_2,f_2)+\zeta^1(a_1,b_1+c_1;d_2,e_2+f_2)+\\
&+\zeta^1(a_1,c_1,b_1;e_2,d_2,f_2)+\zeta^1(a_1+c_1,b_1;d_2+f_2,e_2)+\\
&+\zeta^1(c_1,a_1,b_1;f_2,d_2,e_2)
\end{align*}

Now, let $K$ be a real quadratic field. Let ${\cal{O}}_K$ be the ring of algebraic integers 
in the field $K$.  Let $\mu$ and $\nu$ be totally positive algebraic integers 
in the ring ${\cal{O}}_K$. We define a cone $C$
\[C=\{\alpha\in {\cal{O}}_K\,|\, \alpha=m\mu+n\nu\}.\]
For an element $\alpha$ of the field $K$, 
let $\alpha_1=\sigma_1(\alpha)$ and  $\alpha_2=\sigma_2(\alpha)$ 
be images of the two real embeddings $\sigma_1:K\rightarrow \R$ 
and $\sigma_2:K\rightarrow \R$. 
If $\alpha$ is in $C$ then $\alpha$ is totally positive, 
that is $\alpha_1>0$ and $\alpha_2>0$. 
It is a consequence from the definition of the cone $C$, 
that $\mu$ and $\nu$ are totally positive.
Then we can order elements of the cone $C$ 
in two ways with respect to  each of the two real embeddings. 

Let 
\[\zeta^1_{K,C}(a_1;c_2)
=
\sum_{\alpha\in C}\frac{1}{\alpha_1^{a_1}\alpha_2^{c_2}}
\]
Then we define the following infinite sum shuffle
\begin{align*}
&\zeta^1_{K,C}(a_1;c_2)\zeta^1_{K,C}(b_1;d_2)=\\
=&\sum_{\alpha,\beta\in C}\frac{1}{\alpha_1^{a_1}\beta_1^{b_1}\alpha_2^{c_2}\beta_2^{d_2}}=\\
=&
\left(
\sum_{\mu(\alpha)<\mu(\beta);\,\nu(\alpha)<\nu(\beta)}+
\sum_{\mu(\alpha)=\mu(\beta);\,\nu(\alpha)<\nu(\beta)}+
\sum_{\mu(\alpha)>\mu(\beta);\,\nu(\alpha)<\nu(\beta)}+\right.\\
&+\sum_{\mu(\alpha)<\mu(\beta);\,\nu(\alpha)=\nu(\beta)}+
\sum_{\mu(\alpha)=\mu(\beta);\,\nu(\alpha)=\nu(\beta)}+
\sum_{\mu(\alpha)>\mu(\beta);\,\nu(\alpha)=\nu(\beta)}+\\
&+\left.
\sum_{\mu(\alpha)<\mu(\beta);\,\nu(\alpha)>\nu(\beta)}+
\sum_{\mu(\alpha)=\mu(\beta);\,\nu(\alpha)>\nu(\beta)}+
\sum_{\mu(\alpha)>\mu(\beta);\,\nu(\alpha)>\nu(\beta)}\right)
\frac{1}{\alpha_1^{a_1}\beta_1^{b_1}\alpha_2^{c_2}\beta_2^{d_2}}=\\
=&
\,\,\,_1\zeta_{K,C}(a_1,b_1;c_2,d_2)+\,_{0,1}\zeta_{K,C}(a_1,b_1;c_2,d_2)
+\,_\rho\zeta_{K,C}(a_1,b_1;c_2,d_2)\\
&+\,_{1,0}\zeta_{K,C}(a_1,b_1;c_2,d_2)+\,_1\zeta(a_1+b_1,c_2+d_2)
+\,_{1,0}\zeta_{K,C}(b_1,a_1;d_2,c_2)+\\
&+\,_\rho\zeta_{K,C}(b_1,a_1;d_2,c_2)+\,_{0,1}\zeta_{K,C}(b_1,a_1;d_2,c_2)
+\,_1\zeta_{K,C}(b_1,a_1;d_2,c_2).
\end{align*}

In particular, we obtain
\begin{align}\label{eq stuffle 2,2}
\zeta^1_{K,C}(2;2)\zeta^1_{K,C}(2;2)
=&\zeta^1_{K,C}(4;4)+\,2\,\zeta^1_{K,C}(2,2;2,2)+2\,_\rho\zeta_{K,C}(2,2;2,2)+\\
\nonumber
&+\,_{1,0}\zeta_{K,C}(2,2;2,2)+\,_{1,0}\zeta_{K,C}(2,2;2,2)
\end{align}

\section{Integral shuffle relation}
\subsection{Classical Polylogarithms}
Let us recall the $m$-th polylogarithm and its relation to Riemann zeta values.

The first polylogarithm is defined as 
\[Li_1(x_1)=\int_0^{x_1}\frac{dx_0}{1-x_0}=\int_0^{x_1}(1+x_0+x_0^2+\dots)dx_0=
x_1+\frac{x_1^2}{2}+\frac{x^3_1}{3}+\dots\]
And the $m$-th polylogarithm is defined  by iteration
\begin{equation}
\label{eq n-log}
Li_m(x_m)=\int_0^{x_m}Li_{m-1}(x_{m-1})\frac{dx_{m-1}}{x_{m-1}}.
\end{equation}
Equation \eqref{eq n-log} is a representation of the $m$-th polylogarithm as an iterated integral.
By a direct computation it follows that
$$Li_m(x)=x+\frac{x^2}{2^m}+\frac{x^3}{3^m}+\dots$$
and the relation 
$$\zeta(m)=Li_m(1)$$
is straight forward. 
Using Equation \ref{eq n-log}, we can express the $m$-th polylogarithm as
$$Li_m(x)=\int_{0<x_1<x_2<\dots<x_m<x}
\frac{dx_1}{1-x_1}\wedge\frac{dx_2}{x_2}\wedge\cdots\wedge\frac{dx_{m}}{x_{m}}.$$

Let $x_i=e^{-t_i}$. Then the $m$-th polylogarithm can be written in terms of the variables $t_1,\dots,t_m$ in the following way
\begin{equation}
\label{eq polylog int}
Li_m(e^{-{t}})=\int_{t_1>t_2>\dots>t_m>t}\frac{dt_1\wedge\dots\wedge dt_{m}}{e^{t_1}-1}.
\end{equation}
This can be achieved, first, by changing the variables in the differential forms
$$\frac{dx_1}{1-x_1}=\frac{d(-t_1)}{e^{t_1}-1},\text{ and }\frac{dx_i}{x_i}=d(-t_i),$$
and second, by reversing the bounds of integration $0<x_1<x_2<\dots<x_m<x$ v.s. 
$t_1>t_2>\dots>t_m>t$, which absorbs the sign.
As an infinite sum, we have
\begin{equation}
\label{eq polylog exp}
Li_m(e^{-t})=\sum_{n>0}\frac{e^{-nt}}{n^m}.
\end{equation}

The Equations \eqref{eq polylog int} and  \eqref{eq polylog exp} will be generalized to Dedekind polylogarithms (See Equations \eqref{eq X} and \eqref{eq Y}).

Below we present similar formulas for multiple 
polylogarithms with exponential variables. We will construct  their generalizations to Dedekind multiple polylogarithms (see Lemma \ref{lemma f11}).

Let us recall the definition of double logarithm
\begin{align*}Li_{1,1}(1,x)
&=
\int_0^{x}Li_1(x_1)\frac{dx_1}{1-x_1}
=
\int_0^{x}
\left(\sum_{n_1=1}^{\infty}\frac{x_1^{n_1}}{n_1}\right)
\left(\sum_{n_2=1}^{\infty}x_1^{n_2}\right)
\frac{dx_1}{x_1}
=\\
&=
\sum_{n_1,n_2=1}^{\infty}\frac{x^{n_1+n_2}}{n_1(n_1+n_2)}.
\end{align*}
Let $x_1=e^{-t_1}$ and $x=e^{-t}$. Then the $Li_{1,1}(1,e^{-t})$ can be written as an iterated integral in terms of the variables $t_0,t_1,t$ in the following way:
$$Li_{1,1}(1,e^{-t})=\int_{t_0>t_1>t}\frac{dt_0\wedge dt_1}{(e^{t_0}-1)(e^{t_1}-1)}.$$
As an infinite sum, we have
\begin{equation}
\label{eq double polylog exp}
Li_{1,1}(1,e^{-t})=\sum_{n_1,n_2=1}^\infty\frac{e^{-(n_1+n_2)t}}{n_1(n_1+n_2)}.
\end{equation}

An example of a multiple zeta value is 
$$\zeta(1,2)=\sum_{n_1,n_2=1}^{\infty}\frac{1}{n_1(n_1+n_2)^2}=\int_0^1 Li_{1,1}(x)\frac{dx}{x}.$$

Thus, an integral representation of $\zeta(1,2)$ is
\begin{equation}
\label{eq z12}
\zeta(1,2)
=
\int_{t_1>t_2>t_3>0}\frac{dt_1}{(e^{t_1}-1)}\wedge \frac{dt_2}{(e^{t_2}-1)}\wedge dt_3.
\end{equation}

Similarly,
\begin{equation}
\label{eq z13}
\zeta(1,3)
=
\int_{t_1>t_2>t_3>t_4>0}\frac{dt_1}{(e^{t_1}-1)}\wedge\frac{dt_2}{(e^{t_2}-1)}\wedge dt_3\wedge dt_4.
\end{equation}
and
\begin{equation}
\label{eq z22}
\zeta(2,2)
=
\int_{t_1>t_2>t_3>t_4>0}\frac{dt_1}{(e^{t_1}-1)}\wedge dt_2\wedge\frac{dt_3}{(e^{t_3}-1)}\wedge dt_4.
\end{equation}

Then an integral shuffle relation, sometimes called simply shuffle relation is the following
\begin{align}
\nonumber
\zeta(2)\zeta(2)
=&
\int_{u_1>u_2>0}\frac{du_1}{(e^{u_1}-1)}\wedge du_2
\int_{t_1>t_2>0}\frac{dt_1}{(e^{t_1}-1)}\wedge dt_2=\\
\nonumber
=&
\left(
\int_{u_1>u_2>t_1>t_2>0}+
\int_{u_1>t_1>u_2>t_2>0}\right.+\\
\nonumber
&+\int_{u_1>t_1>t_2>u_2>0}
+\int_{t_1>u_1>u_2>t_2>0}+\\
\nonumber
&+\left.\int_{t_1>u_1>t_2>u_2>0}+
\int_{t_1>t_2>u_1>u_2>0}
\right)
\frac{du_1}{(e^{u_1}-1)}\wedge du_2\wedge\frac{dt_1}{(e^{t_1}-1)}\wedge dt_2=\\
\nonumber
=&\zeta(2,2)+\zeta(1,3)+\zeta(1,3)+\zeta(1,3)+\zeta(1,3)+\zeta(2,2)=\\
\label{eq zeta2,2}=&2\zeta(2,2)+4\zeta(1,3)
\end{align}

From the infinite sum shuffle, sometimes called stuffle relation, we have
\begin{equation}
\label{eq zeta2,2 sum}
\zeta(2)\zeta(2)=\zeta(4)+2\zeta(2,2).
\end{equation}
Therefore, from Equations \eqref{eq zeta2,2} and \eqref{eq zeta2,2 sum}, we obtain the following relation among multiple zeta values.
\begin{equation}
\label{eq relation zeta2,2}
\zeta(1,3)=\frac{1}{4}\zeta(4).
\end{equation}

\subsection{Dedekind Polylogarithms}
A key part of this Subsection is to introduce analogous formulas to \eqref{eq zeta2,2} that provide an integral shuffle for Dedekind zeta values. In order to do that, we need to express (multiple) Dedekind zeta values as  integrals.
First, we are going to recall analogues of (multiple) polylogarithms, which we call Dedekind (multiple) polylogarithms over quadratic extensions on $\Q$ (see Equations \eqref{eq X}, \eqref{eq Y} and Lemma \ref{lemma f11}). 
We will denote by $f_{m}$ the $m$-th Dedekind polylogarithm, which will be an analogue of the $m$-th polylogarithm $Li_m(e^{-t})$ with an exponential variable. Each of the analogues will have an integral representation, resembling an iterated integral and an infinite sum representation, resembling the classical Dedekind zeta values over a quadratic number field. We also draw diagrams that represent  integrals in order to give a geometric view of the iterated integrals on membranes in dimension $2$. We will give examples of multiple Dedekind zeta values (MDZV) over quadratic number field, using the Dedekind (multiple) polylogarithms.

We are going to generalize Equations \eqref{eq polylog exp} and \eqref{eq double polylog exp} for (multiple) polylogarithms to their analogue over the quadratic number field.  We will recall some properties and definitions related to quadratic extensions. For more information one may consider \cite{IR}.

Let $K=\Q\{1,\sqrt{D}\}$. 
By algebraic integers in $K$ we mean all numbers of the form $a+b\sqrt{D}$, where $a$ and $b$ are integers, when $D$ is congruent to 2 or 3 modulo $4$ or numbers of the form
$a+b\frac{1+\sqrt{D}}{2}$. Denote by ${\cal{O}}_K$ the ring of algebraic integers.
 We call the following set $C$ a {\it{cone}}
\[C=\N\{\mu,\nu\}=\{\alpha\in{\cal{O}}_K\mbox{ }|\mbox{ }\alpha=m\mu+n\nu;\mbox{ } m,n \in \N\},\]
where $\N$ denotes the positive integers.
Note that $0$ does not belong to the cone $C$, since the coefficients $a$ and $b$ are positive integers.
If $D>0$ we require that $\mu$ a $\nu$ are totally positive, that is $\mu$, $\nu$ and their Galois conjugates are positive real numbers. If $D<0$, let $\mu$ and $\nu$ have positive real part.
In order to make the examples easier to follow,
we are going to use two sequences of inequalities
$$t_1>u_1>v_1>w_1\mbox{ and }t_2>u_2>v_2>w_2,$$
when we deal with a small number of iterations. 

We are going to define a function $f_1$, which will be an analogue of $Li_1(e^{-t})$. 
Let 
\begin{equation}
\label{eq Gauss f0}
f_0(C;t_1,t_2)=\sum_{\alpha\in C}\exp(-\alpha_1 t_1 - \alpha_2t_2).
\end{equation}
We define $f_1$ as
\begin{equation}
\label{eq f1}
f_1(C,u_1,u_2)=\int^{u_1}_{\infty}\int^{u_2}_{\infty}f_0(C;t_1,t_2)dt_1\wedge dt_2.
\end{equation}
We can draw the following diagram for the integral representing $f_1$.
\begin{center}
\begin{tikzpicture}
\draw[step=2cm] (0,0) grid (2,2);
\draw (0,-.5)node{$+\infty$};
\draw (1,-.5)node{$t_1$};
\draw (2,-.5)node{$u_1$};
\draw (-.5,0)node{$+\infty$};
\draw (-0.5,1)node{$t_2$};
\draw (-.5,2)node{$u_2$};
\draw (1,1)node{$f_0dt_1\wedge dt_2$};
\end{tikzpicture}
\end{center}
The diagram represents that the integrant is $f_0(C;t_1,t_2)dt_1\wedge dt_2$, depending on the variables $t_1$ and $t_2$, subject to the restrictions $+\infty> t_1>u_1$ and $+\infty>t_2>u_2$.

In order to obtain explicit formulas for iterated integrals over membranes, we need the following: 
\begin{lemma}
(a) $$\int^u_{\infty} e^{-kt}dt=\frac{e^{-ku}}{k};$$

(b) Let $N(\alpha)=\alpha_1\alpha_2.$ Then

$$\int^{u_1}_{\infty}\int^{u_2}_{\infty}\exp(-\alpha_1 t_1 - \alpha_2t_2)dt_1\wedge dt_2=
\frac{\exp(-\alpha_1 u_1 - \alpha_2u_2)}{N(\alpha)}.$$
\end{lemma}
The proof is straight forward.

Using the above Lemma, we obtain
$$f_1(C;u_1,u_2)
=
\sum_{\alpha\in C}\frac{\exp(-\alpha_1 u_1- \alpha_2u_2)}{N(\alpha)}
$$

We define a Dedekind dilogarithm $f_2$ by
\begin{align}
\nonumber
f_2(C;v_1,v_2)&=\int^{v_1}_{\infty}\int^{v_2}_{\infty}f_{1}(C;u_1,u_2)d u_1\wedge du_2=\\
\label{eq Ddilog}&=\int_{t_1>u_1>v_1;\mbox{ }t_1>u_1>v_2}f_0(C;t_1,t_2)d t_1\wedge dt_2\wedge d u_1\wedge du_2
\end{align}
We can associate a diagram to the integral representation of the Dedekind dilogarithm $f_2$ (see Equation \eqref{eq Ddilog}).

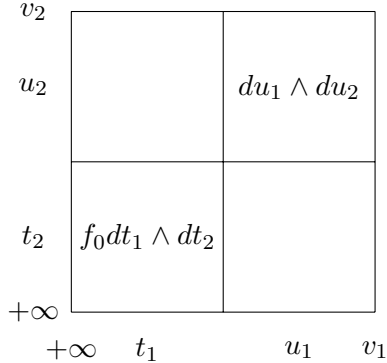
\begin{center}
\begin{tikzpicture}
\draw[step=2cm] (0,0) grid (4,4);
\draw (0,-.5)node{$+\infty$};
\draw (1,-.5)node{$t_1$};
\draw (3,-.5)node{$u_1$};
\draw (4,-.5)node{$v_1$};
\draw (-.5,0)node{$+\infty$};
\draw (-0.5,1)node{$t_2$};
\draw (-.5,3)node{$u_2$};
\draw (-.5,4)node{$v_2$};
\draw (1,1)node{$f_0dt_1\wedge dt_2$};
\draw (3,3)node{$du_1\wedge du_2$};
\end{tikzpicture}
\end{center}
The diagram represents that the variables under the integral are $t_1$, $t_2$, $u_1$, $u_2$, 
subject to the conditions $+\infty>t_1>u_1>v_1$ and $+\infty>t_2>u_2>v_2$. Also, the function $f_0$ in the diagram depends on the variables $t_1$ and $t_2$.

Similarly to Equation \eqref{eq n-log}, we define inductively the $m$-th Dedekind polylogarithm over quadratic number field
\begin{equation}
\label{eq X}
f_m(C;u_1,u_2)=\int^{u_1}_{\infty}\int^{u_1}_{\infty}f_{m-1}(C;t_{1},t_{2})dt_{1}\wedge dt_{2}.
\end{equation}
The above integral is the key example of an iterated integral over a membrane. For more examples and properties one can see \cite{MDZF}.

From Equation \eqref{eq X}, we can derive an analogue of the infinite sum representation of a polylogarithm
(see Equation  \eqref{eq polylog exp}).
\begin{equation}
\label{eq Y}
f_m(C;u_1,u_2)
=
\sum_{\alpha\in C}\frac{\exp(-\alpha u_1 - \overline{\alpha}u_2)}{N(\alpha)^m}.
\end{equation}
The above Equation \eqref{eq Y} gives an infinite sum representation of the $m$-th Dedekind polylogarithm over a quadratic number field.

Now we can define an analogue of the double logarithm $Li_{1,1}(1,e^{-t})$ over the  Gaussian integers, using the following integral representation
$$f_{1,1}(C;v_1,v_2)=\int^{v_1}_{\infty}\int^{v_2}_{\infty}
f_1(C;u_1,u_2)f_0(C;u_1,u_2)du_1\wedge du_2,$$
called a Dedekind double logarithm.
As an analog for Equation \eqref{eq Ddilog}, we can express $f_{1,1}$ only in terms of $f_0$ by
\[f_{1,1}(C;v_1,v_2)=\int_{t_1>u_1>v_1;\mbox{ }t_2>u_2>v_2}(f_0(C;t_1,t_2)d t_1\wedge dt_2)\wedge (f_0(u_1,u_2)d u_1\wedge du_2).\]
It allows us to associate a diagram to the Dedekind double logarithm $f_{1,1}$
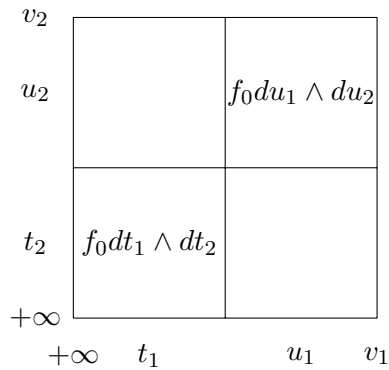
\begin{center}
\begin{tikzpicture}
\draw[step=2cm] (0,0) grid (4,4);
\draw (0,-.5)node{$+\infty$};
\draw (1,-.5)node{$t_1$};
\draw (3,-.5)node{$u_1$};
\draw (4,-.5)node{$v_1$};
\draw (-.5,0)node{$+\infty$};
\draw (-0.5,1)node{$t_2$};
\draw (-.5,3)node{$u_2$};
\draw (-.5,4)node{$v_2$};
\draw (1,1)node{$f_0dt_1\wedge dt_2$};
\draw (3,3)node{$f_0du_1\wedge du_2$};
\end{tikzpicture}
\end{center}
The variables $t_1$, $t_2$, $u_1$, $u_2$ in the diagram are variables in the integrant. They are subject to the conditions $t_1>u_1>v_1$ and $t_2>u_2>v_2$. Also, the lower left function $f_0$ in the diagram depends on the variables $t_1$ and $t_2$, and the upper right function $f_0$ depends on $u_1$ and $u_2$. 

The similarity between $f_{1,1}(C;v_1,v_2)$ and $Li_{1,1}(1,e^{-t_2})$ can be noticed by the infinite sum representation in the following:
\begin{lemma}
\label{lemma f11}
$$f_{1,1}(C;v_1,v_2)
=
\sum_{\alpha,\beta\in C}
\frac
{\exp(-(\alpha+\beta)v_1-(\overline{\alpha}+\overline{\beta})v_2)}
{N(\alpha)N(\alpha+\beta)}.$$
\end{lemma}
\proof
\begin{align*}
f_{1,1}(C;v_1,v_2)&=\int^{v_1}_{\infty}\int^{v_2}_{\infty}
f_1(C;u_1,u_2)f_0(C;u_1,u_2)du_1\wedge du_2=\\
&=\int^{v_1}_{\infty}\int^{v_2}_{\infty}
\sum_{\alpha\in C}
\frac
{\exp(-\alpha_1 u_1-\alpha_2u_2)}
{N(\alpha)}
\sum_{\beta\in C}
\exp(-\beta_1 u_1-\beta_2u_2)
du_1\wedge du_2
=\\
&=\int^{v_1}_{\infty}\int^{v_2}_{\infty}
\sum_{\alpha,\beta\in C}
\frac
{\exp(-(\alpha_1+\beta_1)u_1-(\alpha_2+\beta_2)u_2)}
{N(\alpha)}
du_1\wedge du_2
=\\
&=\sum_{\alpha,\beta\in C}
\frac{\exp(-(\alpha_1+\beta_1)v_1-(\alpha_2+\beta_2)v_2)}
{N(\alpha)N(\alpha+\beta)}.\,\,\,\,\, \mbox{\qed}
\end{align*}

\subsection{Integral shuffle relations}
In this subsection, we consider several examples of shuffle relations based on iterated integrals over membranes.

Let $C_1$ and $C_2$ be two cones. For example,
\[C_1=\N\{1,1+i\}=\{a+b(1+i)\,|\,a,b\in N\}\]
\[C_1=\N\{1,1-i\}=\{a+b(1-i)\,|\,a,b\in N\}\]
Similarly to Equations \eqref{eq Gauss f0} and \eqref{eq f1}, we define 
\[f_0(C_1;t_1,t_2)=\sum_{\alpha\in C_1} \exp(-\alpha_1 t_1 - \alpha_2t_2). \]
and 
\[f_1(C_1;u_1,u_2)=\int^{u_1}_{\infty}\int^{u_2}_{\infty}f_0(C_1;t_1,t_2)dt_1\wedge dt_2.\]

We define Dedekind double logarithms in order to present an integral shuffle relation in a simpler setting. There are two types of Dedekind double logarithms:
\[f_{1,1}(C_1,C_2;v_1,v_2)=\int_{{\bf{t_1>u_1}}>v_1}\int_{t_2>u_2>v_2}f_0(C_1;t_1,t_2)dt_1\wedge dt_2\wedge f_0(C_2;u_1,u_2)du_1\wedge du_2\] and 
\[f_{1,1}^{\rho}(C_1,C_2;v_1,v_2)=
\int_{{\bf{u_1>t_1}}>v_1}\int_{t_2>u_2>v_2}f_0(C_1;t_1,t_2)dt_1\wedge dt_2\wedge f_0(C_2;u_1,t_2)du_1\wedge du_2.\]

\begin{theorem} {\bf{Shuffle relation}} for the Dedekind (poly)-logarithms $f_1(C_1;v_1,v_2)$ and $f_1(C_2;v_1,v_2)$ is the following formula:
\begin{align}
f_1(C_1;v_1,v_2)f_1(C_2;v_1,v_2)
=&f_{1,1}(C_1,C_2;v_1,v_2)+f_{1,1}^{\rho}(C_1,C_2;v_1,v_2)+\\
&+f_{1,1}(C_2,C_1;v_1,v_2)+f_{1,1}^{\rho}(C_2,C_1;v_1,v_2).
\end{align}
\end{theorem}
\proof
\begin{align}
f_1(C_1;v_1,v_2)f_1(C_2;v_1,v_2)
=&\int^{v_1}_{\infty}\int^{v_2}_{\infty}f_0(C;u_1,u_2)du_1\wedge du_2\times\\
&\times\int^{v_1}_{\infty}\int^{v_2}_{\infty}f_0(C';t_1,t_2)dt_1\wedge dt_2=\\
=&\left(\int_{u_1>t_1>v_1}+\int_{t_1>u_1>v_1}\right)
\left(\int_{u_2>t_2>v_2}+\int_{t_2>u_2>v_2}\right)\\
&\,\,f_0(C;u_1,u_2)du_1\wedge du_2\wedge f_0(C';t_1,t_2)dt_1\wedge dt_2=\\
=&f_{1,1}(C,C';v_1,v_2)+f_{1,1}^{\rho}(C,C';v_1,v_2)+\\
&+f_{1,1}(C',C;v_1,v_2)+f_{1,1}^{\rho}(C',C;v_1,v_2).
\end{align}

In terms of diagrams, the shuffle relation can be expressed as
\begin{center}
\begin{tikzpicture}
\draw[step=2cm] (0,0) grid (2,2);
\draw (1,-.5)node{$u_1$};
\draw (2,-.5)node{$v_1$};
\draw (-0.5,1)node{$u_2$};
\draw (-.5,2)node{$v_2$};
\draw (1,1)node{$f_0(C_1)$};
\draw (3,1)node{$\times$};
\draw[step=2cm] (4,0) grid (6,2);
\draw (5,-.5)node{$t_1$};
\draw (6,-.5)node{$v_1$};
\draw (3.5,1)node{$t_2$};
\draw (3.5,2)node{$v_2$};
\draw (5,1)node{$f_0(C_2)$};
\draw (7,1)node{$=$};
\end{tikzpicture}
\end{center}

\begin{center}
\begin{tikzpicture}
\draw[step=2cm] (0,0) grid (4,4);
\draw(-1,2)node{$=$};
\draw (1,-.5)node{$t_1$};
\draw (3,-.5)node{$u_1$};
\draw (4,-.5)node{$v_1$};
\draw (-0.5,1)node{$t_2$};
\draw (-.5,3)node{$u_2$};
\draw (-.5,4)node{$v_2$};
\draw (1,1)node{$f_0(C_1)$};
\draw (3,3)node{$f_0(C_2)$};
\draw(5,2)node{$+$};
\draw[step=2cm] (6,0) grid (10,4);
\draw (7,-.5)node{$t_1$};
\draw (9,-.5)node{$u_1$};
\draw (10,-.5)node{$v_1$};
\draw (5.5,1)node{$t_2$};
\draw (5.5,3)node{$u_2$};
\draw (5.5,4)node{$v_2$};
\draw (9,1)node{$f_0(C_1)$};
\draw (7,3)node{$f_0(C_2)$};
\draw (11,2)node{$+$};
\end{tikzpicture}
\end{center}

\begin{center}
\begin{tikzpicture}
\draw[step=2cm] (0,0) grid (4,4);
\draw(-1,2)node{$+$};
\draw (1,-.5)node{$t_1$};
\draw (3,-.5)node{$u_1$};
\draw (4,-.5)node{$v_1$};
\draw (-0.5,1)node{$t_2$};
\draw (-.5,3)node{$u_2$};
\draw (-.5,4)node{$v_2$};
\draw (1,1)node{$f_0(C_2)$};
\draw (3,3)node{$f_0(C_1)$};
\draw(5,2)node{$+$};
\draw[step=2cm] (6,0) grid (10,4);
\draw (7,-.5)node{$t_1$};
\draw (9,-.5)node{$u_1$};
\draw (10,-.5)node{$v_1$};
\draw (5.5,1)node{$t_2$};
\draw (5.5,3)node{$u_2$};
\draw (5.5,4)node{$v_2$};
\draw (9,1)node{$f_0(C_2)$};
\draw (7,3)node{$f_0(C_1)$};
\end{tikzpicture}
\end{center}

\qed

Similarly to the Dedekind double logarithm $f_{1,1}$, we define
$$f_{1,2}(C,w_1,w_2)
=
\int^{w_1}_{\infty}\int^{w_2}_{\infty}
f_{1,1}(C;v_1,v_2)dv_1\wedge dv_2.$$
We can associate the following diagram to the multiple Dedekind polylogarithm $f_{1,2}$

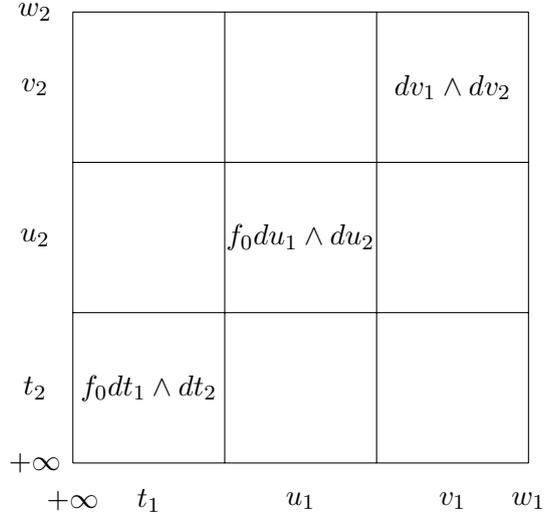
\begin{center}
\begin{tikzpicture}
\draw[step=2cm] (0,0) grid (6,6);
\draw (0,-.5)node{$+\infty$};
\draw (1,-.5)node{$t_1$};
\draw (3,-.5)node{$u_1$};
\draw (5,-.5)node{$v_1$};
\draw (6,-.5)node{$w_1$};
\draw (-0.5,0)node{$+\infty$};
\draw (-0.5,1)node{$t_2$};
\draw (-.5,3)node{$u_2$};
\draw (-.5,5)node{$v_2$};
\draw (-.5,6)node{$w_2$};
\draw (1,1)node{$f_0dt_1\wedge dt_2$};
\draw (3,3)node{$f_0du_1\wedge du_2$};
\draw (5,5)node{$dv_1\wedge dv_2$};
\end{tikzpicture}
\end{center}
The diagram represents the following: The variables of the integrant are $t_1$,  $t_2$, $u_1$, $u_2$, $v_1$,  $v_2$. 
The variables are subject to the conditions $t_1>u_1>v_1>w_1$ and $t_2>u_2>v_2>w_2$.
The lower left function $f_0$ depends on the variables $t_1$ and $t_2$. And the middle function $f_0$ depends on $u_1$ and $u_2$. The upper right $2$-form is $dv_1\wedge dv_2$. Thus, the diagram represents the following integral:
\begin{align}
\label{eq pre f12}
&f_{1,2}(C;w_1,w_2)=\\
\nonumber
&=\int_{D_{w_1,w_2}}
(f_0(C;t_1,t_2)dt_1\wedge dt_2)\wedge
(f_0(C;u_1,u_2)du_1\wedge du_2)\wedge (dv_1\wedge dv_2),
\end{align}
where the domain of integration is
\[D_{w_1,w_2}=\{(t_1,t_2,u_1,u_2,v_1,v_2)\in \R^6\mbox{ }|\mbox{ }t_1>u_1>v_1>w_1\mbox{ and } t_2>u_2>v_2>w_2\}\]
A direct computation leads to
$$f_{1,2}(C;w_1,w_2)
=\sum_{\alpha,\beta\in C}
\frac{\exp(-(\alpha+\beta)w_1-(\overline{\alpha}+\overline{\beta})w_2)}
{N(\alpha)N(\alpha+\beta)^2}.
$$
Using the above Equation, we define a multiple Dedekind zeta value as 
$$\zeta_{\Q(i);C}(1,2;1,2)
=
f_{1,2}(C;0,0)
=
\sum_{\alpha,\beta\in C}
\frac{1}{N(\alpha)N(\alpha+\beta)^2}.
$$

Similarly to the Dedekind polylogarithm $f_{1,2}$, we define
$$f_{1,2}^{(1,2)(1)}(C,w_1,w_2)
=
\int^{w_1}_{\infty}\int^{w_2}_{\infty}
f_{1,1}^{(1,2),(1)}(C,C;v_1,v_2)dv_1\wedge dv_2.$$
The definition above is needed for the integral shuffle relation for multiple Dedekind zeta values. We can associate the following diagram to the multiple Dedekind polylogarithm $f_{1,2}^{(1,2),(1)}$

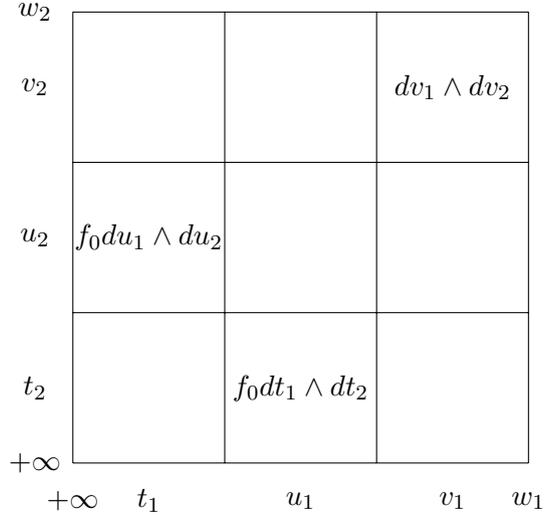
\begin{center}
\begin{tikzpicture}
\draw[step=2cm] (0,0) grid (6,6);
\draw (0,-.5)node{$+\infty$};
\draw (1,-.5)node{$t_1$};
\draw (3,-.5)node{$u_1$};
\draw (5,-.5)node{$v_1$};
\draw (6,-.5)node{$w_1$};
\draw (-0.5,0)node{$+\infty$};
\draw (-0.5,1)node{$t_2$};
\draw (-.5,3)node{$u_2$};
\draw (-.5,5)node{$v_2$};
\draw (-.5,6)node{$w_2$};
\draw (3,1)node{$f_0dt_1\wedge dt_2$};
\draw (1,3)node{$f_0du_1\wedge du_2$};
\draw (5,5)node{$dv_1\wedge dv_2$};
\end{tikzpicture}
\end{center}
The diagram represents the following: The variables of the integrant are $t_1$,  $t_2$, $u_1$, $u_2$, $v_1$,  $v_2$. 
The variables are subject to the conditions $t_1>u_1>v_1>w_1$ and $t_2>u_2>v_2>w_2$.
One of the functions $f_0$ depends on the variables $u_1$ and $t_2$. And the other function $f_0$ depends on $t_1$ and $u_2$.  

Thus, from the diagram  we obtain the following integral, which is more useful for writing explicit formulas:
\begin{align}
\label{eq pre f12}
&f_{1,2}^{(1,2),(1)}(C;w_1,w_2)=\\
\nonumber
&=\int_{D_{w_1,w_2}}
(f_0(C;t_1,t_2)dt_1\wedge dt_2)\wedge
(f_0(C;u_1,u_2)du_1\wedge du_2)\wedge (dv_1\wedge dv_2),
\end{align}
where the domain of integration is
\[D_{w_1,w_2}=\{(t_1,t_2,u_1,u_2,v_1,v_2)\in \R^6\mbox{ }|\mbox{ }{\bf{u_1>t_1}}>v_1>w_1\mbox{ and } t_2>u_2>v_2>w_2\}\]
A direct computation leads to
$$f_{1,2}^{(1,2),(1)}(C;w_1,w_2)
=\sum_{\alpha,\beta\in C}
\frac{\exp(-(\alpha+\beta)w_1-(\overline{\alpha}+\overline{\beta})w_2)}
{\alpha_1\beta_2N(\alpha+\beta)^2}.
$$
We define a multiple Dedekind zeta value as 
$$\zeta_{\Q(i);C}^{(1,2),(1)}(1,2;1,2)
=
f_{1,2}^{(1,2),(1)}(C;0,0)
=
\sum_{\alpha,\beta\in C}
\frac{1}{\beta_1(\alpha_1+\beta_1)^2\alpha_2(\alpha_2+\beta_2)^2}.
$$
On pages 20-23, we will consider 36 similar diagrams.  Then, the above diagram associated to 
$\zeta_{\Q(i);C}^{(1,2),(1)}(1,2;1,2)$ would be denoted by
\begin{center}
\begin{tikzpicture}
\draw[step=1cm] (0,0) grid (3,3); 
\draw (.5,1.5) node{$1$};
\draw (1.5,.5) node{$1$};
\draw (2.5,2.5)node{$0$};
\draw (.5, -.5)node{$t_1$};
\draw (1.5, -.5)node{$u_1$};
\draw (2.5, -.5)node{$v_1$};
\draw (-.5, .5)node{$u_2$};
\draw (-.5, 1.5)node{$t_2$};
\draw (-.5, 2.5)node{$v_2$};
\end{tikzpicture}
\end{center}
The $1$'s signify $f_0(t_1,t_2)dt_1\wedge dt_2$ or $f_0(u_1,u_2)du_1\wedge du_2$. The $0$ signifies the differential form $dv_1\wedge dv_2$. Also, for the valuables $t_1,u_1,v_1$ we have $t_1>u_1v>1>0$. Similarly, $u_2>t_2>v_2>0$ is the meaning of the variables left of the box.

Compare the above diagram to the previous one. They denote the same integrals, when $w_1=0$ and $w_2=0$

Now, we are going to consider the particular example of an integral shuffle relation for the product
$\zeta_{K,C}(2)\zeta_{K,C}(2).$

Let ${\cal{O}}_K$ be the ring of integers in a real quadratic field $K$.
Let $C$ be a cone
\[C=\N\{\mu,\nu\}=\{m\mu+n\nu\,|\,m,n\in \N\},\]
where both $\mu$ and $\nu$ are totally positive integers.

Then, we have the following integral shuffle:

\begin{theorem} (Integral shuffle)
\label{theorem int shuffle 2,2}
\begin{align}
\label{eq int shuffle 2,2}
\zeta_{K,C}(2)\zeta_{K,C}(2)
&=2\zeta^1_{K,C}(2,2;2,2)+8\zeta^1_{K,C}(1,3;1,3)+\\
\nonumber
&+4\zeta^{1}_{K,C}(1,3;2,2)+4\zeta^{1}_{K,C}(2,2;1,3)+\\
\nonumber
&+2\zeta^\rho_{K,C}(2,2;2,2)+8\zeta^\rho_{K,C}(1,3;1,3)+\\
\nonumber
&+4\zeta^\rho_{K,C}(1,3;2,2)+4\zeta^\rho_{K,C}(2,2;1,3).
\end{align}
\end{theorem}
\proof
The shuffle is done in the following way. The variables $t_1>u_1>0$ are shuffled with the variables $v_1>w_1>0$. There are $6$ such possibilities. Then variables $t_2>u_2>0$ are shuffled with the variables $v_2>w_2>0$. Shuffling simultaneously among the variables $t_1,u_1,v_1,w_1$ with index $1$ and among the variables  $t_2,u_2,v_2,w_2$  with index $2$, we obtain $6^2$ possibilities. When the variables with index $1$ are shuffled, we have $2$ summation of  the type $\frac{1}{\alpha_1^2(\alpha_1+\beta_1)^2}\times \cdots$ and $4$ summations of the type $\frac{1}{\alpha_1(\alpha_1+\beta_1)^3}\times \cdots$. We obtain similar terms when we shuffle the variables with index $2$, namely,  $2$ summation of  the type $\frac{1}{\alpha_2^2(\alpha_2+\beta_2)^2}\times \cdots$ and $4$ summations of the type $\frac{1}{\alpha_2(\alpha_2+\beta_2)^3}\times \cdots$.

\begin{align}
\nonumber
\zeta_{K,C}(2)\zeta_{K,C}(2)=&\int_{t_1>u_1>0;\,\,t_2>u_2>0}f_0(C;t_1,t_2)dt_1\wedge dt_2\wedge du_1\wedge du_2\times\\
\nonumber
&\times\int_{v_1>w_1>0;\,\,v_2>w_2>0}f_0(C;v_1,v_2)dv_1\wedge dv_2\wedge dw_1\wedge dw_2=\\
\label{eq int shuffle 2,2}
=&2\sum_{\alpha,\beta\in C}\frac{1}{N(\alpha)^2N(\alpha+\beta)^2}
+8\sum_{\alpha,\beta\in C}\frac{1}{N(\alpha)N(\alpha+\beta)^3}+\\
\nonumber
&+4\sum_{\alpha,\beta\in C}
\frac{1}{\alpha_1^1(\alpha_1+\beta_1)^3\alpha_2^2(\alpha_2+\beta_2)^2}+\\
\nonumber
&+4\sum_{\alpha,\beta\in C}
\frac{1}{\alpha_1^2(\alpha_1+\beta_1)^2\alpha_2^1(\alpha_2+\beta_2)^3}+\\
\nonumber
&+2\sum_{\alpha,\beta\in C}\frac{1}{\alpha_1^2\beta_2^2N(\alpha+\beta)^2}
+8\sum_{\alpha,\beta\in C}\frac{1}{\alpha_1\beta_2N(\alpha+\beta)^3}+\\
\nonumber
&+4\sum_{\alpha,\beta\in C}\frac{1}{\alpha_1^2(\alpha_1+\beta_1)^2\beta_2^1(\alpha_2+\beta_2)^3}+\\
\nonumber
&+4\sum_{\alpha,\beta\in C}\frac{1}{\alpha_1^1(\alpha_1+\beta_1)^3\beta_2^2(\alpha_2+\beta_2)^2}=\\
=&2\zeta^1_{K,C}(2,2;2,2)+8\zeta^1_{K,C}(1,3;1,3)+\\
&+4\zeta^{1}_{K,C}(1,3;2,2)+4\zeta^{1}_{K,C}(2,2;1,3)+\\
&+2\zeta^\rho_{K,C}(2,2;2,2)+8\zeta^\rho_{K,C}(1,3;1,3)+\\
&+4\zeta^\rho_{K,C}(1,3;2,2)+4\zeta^\rho_{K,C}(2,2;1,3).\\
\end{align}

The following 36 diagrams reprepresent all the possible shuffles among the variables $t_1>u_1>0$
and $v_1>w_1>0$, and among $t_2>u_2>0$ and $v_2>w_2>0$. The diagrams below show all possibilities of how these variables can be ordered respecting the given inequalities: $t_1>u_1>0$, $v_1>w_1>0$, $t_2>u_2>0$ and $v_2>w_2>0$. In the boxes, the number $1$ appears at coordinate $(t_1,t_2)$ and at $(v_1,v_2)$, they correspond to the 2-forms $f_0(t_1,t_2)dt_1\wedge dt_2$ and $f_0(v_1,v_2)dv_1\wedge dv_2$ under the integral. Also, the occurrence of $0$ in the boxes appears at coordinates $(u_1,u_2)$ and $(w_1,w_2)$. It signifies that under the integral we have the forms $du_1\wedge du_2$ and $dw_1\wedge dw_2$, respectively. In each diagram, the order of the valuables in horizontal and vertical direction are in decreasing order of the variables.

One term of the infinite sum is given below each of the diagrams. It is one term of a multiple Dedekind zeta value from the infinite sum representation. Note that a diagram itself denotes an integral representation of a multiple Dedekind zeta value. 



\begin{center}
\begin{tikzpicture}
\draw[step=1cm] (-5,0) grid (-1,4);
\draw[step=1cm] (0,0) grid (4,4);
\draw[step=1cm] (5,0) grid (9,4);
\draw (-4.5,-.5)node{$t_1$};
\draw (-3.5,-.5)node{$u_1$};
\draw (-2.5,-.5)node{$v_1$};
\draw (-1.5,-.5)node{$w_1$};

\draw (.5,-.5)node{$t_1$};
\draw (1.5,-.5)node{$v_1$};
\draw (2.5,-.5)node{$u_1$};
\draw (3.5,-.5)node{$w_1$};

\draw (5.5,-.5)node{$t_1$};
\draw (6.5,-.5)node{$v_1$};
\draw (7.5,-.5)node{$w_1$};
\draw (8.5,-.5)node{$u_1$};
\draw (-5.3,.5)node{$t_2$};
\draw (-5.3,1.5)node{$u_2$};
\draw (-5.3,2.5)node{$v_2$};
\draw (-5.3,3.5)node{$w_2$};

\draw (-.3,.5)node{$t_2$};
\draw (-.3,1.5)node{$u_2$};
\draw (-.3,2.5)node{$v_2$};
\draw (-.3,3.5)node{$w_2$};

\draw (4.7,.5)node{$t_2$};
\draw (4.7,1.5)node{$u_2$};
\draw (4.7,2.5)node{$v_2$};
\draw (4.7,3.5)node{$w_2$};
\draw (-4.5,.5)node{$1$};
\draw (-3.5,1.5)node{$0$};
\draw (-2.5,2.5)node{$1$};
\draw (-1.5,3.5)node{$0$};

\draw (.5,.5)node{$1$};
\draw (2.5,1.5)node{$0$};
\draw (1.5,2.5)node{$1$};
\draw (3.5,3.5)node{$0$};

\draw (5.5,.5)node{$1$};
\draw (7.5,3.5)node{$0$};
\draw (6.5,2.5)node{$1$};
\draw (8.5,1.5)node{$0$};
\draw (-3,-1.5)node{$\frac{1}{N(\alpha)^2N(\alpha+\beta)^2}$};
\draw (2,-1.5)node{$\frac{1}{\alpha_1(\alpha_1+\beta_1)^3\alpha_2^2(\alpha_2+\beta_2)^2}$};
\draw (7,-1.5)node{$\frac{1}{\alpha_1(\alpha_1+\beta_1)^3\alpha_2^2(\alpha_2+\beta_2)^2}$};
\draw[step=1cm] (-5,-7) grid (-1,-3);
\draw[step=1cm] (0,-7) grid (4,-3);
\draw[step=1cm] (5,-7) grid (9,-3);
\draw (-4.5,-7.5)node{$t_1$};
\draw (-3.5,-7.5)node{$u_1$};
\draw (-2.5,-7.5)node{$v_1$};
\draw (-1.5,-7.5)node{$w_1$};

\draw (.5,-7.5)node{$t_1$};
\draw (1.5,-7.5)node{$v_1$};
\draw (2.5,-7.5)node{$u_1$};
\draw (3.5,-7.5)node{$w_1$};

\draw (5.5,-7.5)node{$t_1$};
\draw (6.5,-7.5)node{$v_1$};
\draw (7.5,-7.5)node{$w_1$};
\draw (8.5,-7.5)node{$u_1$};
\draw (-5.3,-6.5)node{$t_2$};
\draw (-5.3,-5.5)node{$v_2$};
\draw (-5.3,-4.5)node{$u_2$};
\draw (-5.3,-3.5)node{$w_2$};

\draw (-.3,-6.5)node{$t_2$};
\draw (-.3,-5.5)node{$v_2$};
\draw (-.3,-4.5)node{$u_2$};
\draw (-.3,-3.5)node{$w_2$};

\draw (4.7,-6.5)node{$t_2$};
\draw (4.7,-5.5)node{$v_2$};
\draw (4.7,-4.5)node{$u_2$};
\draw (4.7,-3.5)node{$w_2$};
\draw (-4.5,-6.5)node{$1$};
\draw (-2.5,-5.5)node{$1$};
\draw (-3.5,-4.5)node{$0$};
\draw (-1.5,-3.5)node{$0$};

\draw (.5,-6.5)node{$1$};
\draw (1.5,-5.5)node{$1$};
\draw (2.5,-4.5)node{$0$};
\draw (3.5,-3.5)node{$0$};

\draw (5.5,-6.5)node{$1$};
\draw (6.5,-5.5)node{$1$};
\draw (7.5,-3.5)node{$0$};
\draw (8.5,-4.5)node{$0$};
\draw (-3,-8.5)node{$\frac{1}{\alpha_1^2(\alpha_1+\beta_1)^2\alpha_2(\alpha_2+\beta_2)^3}$};
\draw (2,-8.5)node{$\frac{1}{N(\alpha)^1N(\alpha+\beta)^3}$};
\draw (7,-8.5)node{$\frac{1}{N(\alpha)^1N(\alpha+\beta)^3}$};

\draw[step=1cm] (-5,-14) grid (-1,-10);
\draw[step=1cm] (0,-14) grid (4,-10);
\draw[step=1cm] (5,-14) grid (9,-10);
\draw (-4.5,-14.5)node{$t_1$};
\draw (-3.5,-14.5)node{$u_1$};
\draw (-2.5,-14.5)node{$v_1$};
\draw (-1.5,-14.5)node{$w_1$};

\draw (.5,-14.5)node{$t_1$};
\draw (1.5,-14.5)node{$v_1$};
\draw (2.5,-14.5)node{$u_1$};
\draw (3.5,-14.5)node{$w_1$};

\draw (5.5,-14.5)node{$t_1$};
\draw (6.5,-14.5)node{$v_1$};
\draw (7.5,-14.5)node{$w_1$};
\draw (8.5,-14.5)node{$u_1$};
\draw (-5.3,-13.5)node{$t_2$};
\draw (-5.3,-12.5)node{$v_2$};
\draw (-5.3,-11.5)node{$w_2$};
\draw (-5.3,-10.5)node{$u_2$};

\draw (-.3,-13.5)node{$t_2$};
\draw (-.3,-12.5)node{$v_2$};
\draw (-.3,-11.5)node{$w_2$};
\draw (-.3,-10.5)node{$u_2$};

\draw (4.7,-13.5)node{$t_2$};
\draw (4.7,-12.5)node{$v_2$};
\draw (4.7,-11.5)node{$w_2$};
\draw (4.7,-10.5)node{$u_2$};
\draw (-4.5,-13.5)node{$1$};
\draw (-2.5,-12.5)node{$1$};
\draw (-1.5,-11.5)node{$0$};
\draw (-3.5,-10.5)node{$0$};

\draw (.5,-13.5)node{$1$};
\draw (1.5,-12.5)node{$1$};
\draw (3.5,-11.5)node{$0$};
\draw (2.5,-10.5)node{$0$};

\draw (5.5,-13.5)node{$1$};
\draw (6.5,-12.5)node{$1$};
\draw (7.5,-11.5)node{$0$};
\draw (8.5,-10.5)node{$0$};
\draw (-3,-15.5)node{$\frac{1}{\alpha_1^2(\alpha_1+\beta_1)^2\alpha_2(\alpha_2+\beta_2)^3}$};
\draw (2,-15.5)node{$\frac{1}{N(\alpha)^1N(\alpha+\beta)^3}$};
\draw (7,-15.5)node{$\frac{1}{N(\alpha)^1N(\alpha+\beta)^3}$};
\draw(2,-16.5)node{{\bf{Table (1,1)}}};
\end{tikzpicture}

\end{center}

\begin{center}
\begin{tikzpicture}
\draw[step=1cm] (-5,0) grid (-1,4);
\draw[step=1cm] (0,0) grid (4,4);
\draw[step=1cm] (5,0) grid (9,4);
\draw (-4.5,-.5)node{$v_1$};
\draw (-3.5,-.5)node{$t_1$};
\draw (-2.5,-.5)node{$u_1$};
\draw (-1.5,-.5)node{$w_1$};

\draw (.5,-.5)node{$v_1$};
\draw (1.5,-.5)node{$t_1$};
\draw (2.5,-.5)node{$w_1$};
\draw (3.5,-.5)node{$u_1$};

\draw (5.5,-.5)node{$v_1$};
\draw (6.5,-.5)node{$w_1$};
\draw (7.5,-.5)node{$t_1$};
\draw (8.5,-.5)node{$u_1$};
\draw (-5.3,.5)node{$t_2$};
\draw (-5.3,1.5)node{$u_2$};
\draw (-5.3,2.5)node{$v_2$};
\draw (-5.3,3.5)node{$w_2$};

\draw (-.3,.5)node{$t_2$};
\draw (-.3,1.5)node{$u_2$};
\draw (-.3,2.5)node{$v_2$};
\draw (-.3,3.5)node{$w_2$};

\draw (4.7,.5)node{$t_2$};
\draw (4.7,1.5)node{$u_2$};
\draw (4.7,2.5)node{$v_2$};
\draw (4.7,3.5)node{$w_2$};
\draw (-3.5,.5)node{$1$};
\draw (-2.5,1.5)node{$0$};
\draw (-4.5,2.5)node{$1$};
\draw (-1.5,3.5)node{$0$};

\draw (1.5,.5)node{$1$};
\draw (2.5,3.5)node{$0$};
\draw (.5,2.5)node{$1$};
\draw (3.5,1.5)node{$0$};

\draw (7.5,.5)node{$1$};
\draw (6.5,3.5)node{$0$};
\draw (5.5,2.5)node{$1$};
\draw (8.5,1.5)node{$0$};
\draw (-3,-1.5)node{$\frac{1}{\beta_1(\alpha_1+\beta_1)^3\alpha_2^2(\alpha_2+\beta_2)^2}$};
\draw (2,-1.5)node{$\frac{1}{\beta_1(\alpha_1+\beta_1)^3\alpha_2^2(\alpha_2+\beta_2)^2}$};
\draw (7,-1.5)node{$\frac{1}{\beta_1^2(\alpha_1+\beta_1)^2\alpha_2^2(\alpha_2+\beta_2)^2}$};
\draw[step=1cm] (-5,-7) grid (-1,-3);
\draw[step=1cm] (0,-7) grid (4,-3);
\draw[step=1cm] (5,-7) grid (9,-3);
\draw (-4.5,-7.5)node{$v_1$};
\draw (-3.5,-7.5)node{$t_1$};
\draw (-2.5,-7.5)node{$u_1$};
\draw (-1.5,-7.5)node{$w_1$};

\draw (.5,-7.5)node{$v_1$};
\draw (1.5,-7.5)node{$t_1$};
\draw (2.5,-7.5)node{$w_1$};
\draw (3.5,-7.5)node{$u_1$};

\draw (5.5,-7.5)node{$v_1$};
\draw (6.5,-7.5)node{$w_1$};
\draw (7.5,-7.5)node{$t_1$};
\draw (8.5,-7.5)node{$u_1$};
\draw (-5.3,-6.5)node{$t_2$};
\draw (-5.3,-5.5)node{$v_2$};
\draw (-5.3,-4.5)node{$u_2$};
\draw (-5.3,-3.5)node{$w_2$};

\draw (-.3,-6.5)node{$t_2$};
\draw (-.3,-5.5)node{$v_2$};
\draw (-.3,-4.5)node{$u_2$};
\draw (-.3,-3.5)node{$w_2$};

\draw (4.7,-6.5)node{$t_2$};
\draw (4.7,-5.5)node{$v_2$};
\draw (4.7,-4.5)node{$u_2$};
\draw (4.7,-3.5)node{$w_2$};
\draw (-3.5,-6.5)node{$1$};
\draw (-4.5,-5.5)node{$1$};
\draw (-2.5,-4.5)node{$0$};
\draw (-1.5,-3.5)node{$0$};

\draw (1.5,-6.5)node{$1$};
\draw (.5,-5.5)node{$1$};
\draw (3.5,-4.5)node{$0$};
\draw (2.5,-3.5)node{$0$};

\draw (7.5,-6.5)node{$1$};
\draw (5.5,-5.5)node{$1$};
\draw (6.5,-3.5)node{$0$};
\draw (8.5,-4.5)node{$0$};
\draw (-3,-8.5)node{$\frac{1}{\beta_1(\alpha_1+\beta_1)^3\alpha_2(\alpha_2+\beta_2)^3}$};
\draw (2,-8.5)node{$\frac{1}{\beta_1(\alpha_1+\beta_1)^3\alpha_2(\alpha_2+\beta_2)^3}$};
\draw (7,-8.5)node{$\frac{1}{\beta_1^2(\alpha_1+\beta_1)^2\alpha_2^1(\alpha_2+\beta_2)^3}$};

\draw[step=1cm] (-5,-14) grid (-1,-10);
\draw[step=1cm] (0,-14) grid (4,-10);
\draw[step=1cm] (5,-14) grid (9,-10);
\draw (-4.5,-14.5)node{$v_1$};
\draw (-3.5,-14.5)node{$t_1$};
\draw (-2.5,-14.5)node{$u_1$};
\draw (-1.5,-14.5)node{$w_1$};

\draw (.5,-14.5)node{$v_1$};
\draw (1.5,-14.5)node{$t_1$};
\draw (2.5,-14.5)node{$w_1$};
\draw (3.5,-14.5)node{$u_1$};

\draw (5.5,-14.5)node{$v_1$};
\draw (6.5,-14.5)node{$w_1$};
\draw (7.5,-14.5)node{$t_1$};
\draw (8.5,-14.5)node{$u_1$};
\draw (-5.3,-13.5)node{$t_2$};
\draw (-5.3,-12.5)node{$v_2$};
\draw (-5.3,-11.5)node{$w_2$};
\draw (-5.3,-10.5)node{$u_2$};

\draw (-.3,-13.5)node{$t_2$};
\draw (-.3,-12.5)node{$v_2$};
\draw (-.3,-11.5)node{$w_2$};
\draw (-.3,-10.5)node{$u_2$};

\draw (4.7,-13.5)node{$t_2$};
\draw (4.7,-12.5)node{$v_2$};
\draw (4.7,-11.5)node{$w_2$};
\draw (4.7,-10.5)node{$u_2$};
\draw (-3.5,-13.5)node{$1$};
\draw (-4.5,-12.5)node{$1$};
\draw (-1.5,-11.5)node{$0$};
\draw (-2.5,-10.5)node{$0$};

\draw (1.5,-13.5)node{$1$};
\draw (.5,-12.5)node{$1$};
\draw (2.5,-11.5)node{$0$};
\draw (3.5,-10.5)node{$0$};

\draw (7.5,-13.5)node{$1$};
\draw (5.5,-12.5)node{$1$};
\draw (6.5,-11.5)node{$0$};
\draw (8.5,-10.5)node{$0$};
\draw (-3,-15.5)node{$\frac{1}{\beta_1(\alpha_1+\beta_1)^3\alpha_2(\alpha_2+\beta_2)^3}$};
\draw (2,-15.5)node{$\frac{1}{\beta_1(\alpha_1+\beta_1)^3\alpha_2(\alpha_2+\beta_2)^3}$};
\draw (7,-15.5)node{$\frac{1}{\beta_1^2(\alpha_1+\beta_1)^2\alpha_2^1(\alpha_2+\beta_2)^3}$};
\draw(2,-16.5)node{{\bf{Table (1,2)}}};
\end{tikzpicture}
\end{center}


\begin{center}
\begin{tikzpicture}
\draw[step=1cm] (-5,0) grid (-1,4);
\draw[step=1cm] (0,0) grid (4,4);
\draw[step=1cm] (5,0) grid (9,4);
\draw (-4.5,-.5)node{$t_1$};
\draw (-3.5,-.5)node{$u_1$};
\draw (-2.5,-.5)node{$v_1$};
\draw (-1.5,-.5)node{$w_1$};

\draw (.5,-.5)node{$t_1$};
\draw (1.5,-.5)node{$v_1$};
\draw (2.5,-.5)node{$u_1$};
\draw (3.5,-.5)node{$w_1$};

\draw (5.5,-.5)node{$t_1$};
\draw (6.5,-.5)node{$v_1$};
\draw (7.5,-.5)node{$w_1$};
\draw (8.5,-.5)node{$u_1$};
\draw (-5.3,.5)node{$v_2$};
\draw (-5.3,1.5)node{$t_2$};
\draw (-5.3,2.5)node{$u_2$};
\draw (-5.3,3.5)node{$w_2$};

\draw (-.3,.5)node{$v_2$};
\draw (-.3,1.5)node{$t_2$};
\draw (-.3,2.5)node{$u_2$};
\draw (-.3,3.5)node{$w_2$};

\draw (4.7,.5)node{$v_2$};
\draw (4.7,1.5)node{$t_2$};
\draw (4.7,2.5)node{$u_2$};
\draw (4.7,3.5)node{$w_2$};
\draw (-2.5,.5)node{$1$};
\draw (-4.5,1.5)node{$1$};
\draw (-3.5,2.5)node{$0$};
\draw (-1.5,3.5)node{$0$};

\draw (1.5,.5)node{$1$};
\draw (.5,1.5)node{$1$};
\draw (2.5,2.5)node{$0$};
\draw (3.5,3.5)node{$0$};

\draw (6.5,.5)node{$1$};
\draw (7.5,3.5)node{$0$};
\draw (5.5,1.5)node{$1$};
\draw (8.5,2.5)node{$0$};
\draw (-3,-1.5)node{$\frac{1}{\beta_1^2(\alpha_1+\beta_1)^2\alpha_2(\alpha_2+\beta_2)^3}$};
\draw (2,-1.5)node{$\frac{1}{\beta_1(\alpha_1+\beta_1)^3\alpha_2(\alpha_2+\beta_2)^3}$};
\draw (7,-1.5)node{$\frac{1}{\beta_1(\alpha_1+\beta_1)^3\alpha_2(\alpha_2+\beta_2)^3}$};
\draw[step=1cm] (-5,-7) grid (-1,-3);
\draw[step=1cm] (0,-7) grid (4,-3);
\draw[step=1cm] (5,-7) grid (9,-3);
\draw (-4.5,-7.5)node{$t_1$};
\draw (-3.5,-7.5)node{$u_1$};
\draw (-2.5,-7.5)node{$v_1$};
\draw (-1.5,-7.5)node{$w_1$};

\draw (.5,-7.5)node{$t_1$};
\draw (1.5,-7.5)node{$v_1$};
\draw (2.5,-7.5)node{$u_1$};
\draw (3.5,-7.5)node{$w_1$};

\draw (5.5,-7.5)node{$t_1$};
\draw (6.5,-7.5)node{$v_1$};
\draw (7.5,-7.5)node{$w_1$};
\draw (8.5,-7.5)node{$u_1$};
\draw (-5.3,-6.5)node{$v_2$};
\draw (-5.3,-5.5)node{$t_2$};
\draw (-5.3,-4.5)node{$w_2$};
\draw (-5.3,-3.5)node{$u_2$};

\draw (-.3,-6.5)node{$v_2$};
\draw (-.3,-5.5)node{$t_2$};
\draw (-.3,-4.5)node{$w_2$};
\draw (-.3,-3.5)node{$u_2$};

\draw (4.7,-6.5)node{$v_2$};
\draw (4.7,-5.5)node{$t_2$};
\draw (4.7,-4.5)node{$w_2$};
\draw (4.7,-3.5)node{$u_2$};
\draw (-2.5,-6.5)node{$1$};
\draw (-4.5,-5.5)node{$1$};
\draw (-1.5,-4.5)node{$0$};
\draw (-3.5,-3.5)node{$0$};

\draw (1.5,-6.5)node{$1$};
\draw (.5,-5.5)node{$1$};
\draw (3.5,-4.5)node{$0$};
\draw (2.5,-3.5)node{$0$};

\draw (6.5,-6.5)node{$1$};
\draw (5.5,-5.5)node{$1$};
\draw (8.5,-3.5)node{$0$};
\draw (7.5,-4.5)node{$0$};
\draw (-3,-8.5)node{$\frac{1}{\alpha_1^2(\alpha_1+\beta_1)^2\beta_2(\alpha_2+\beta_2)^3}$};
\draw (2,-8.5)node{$\frac{1}{\beta_1(\alpha_1+\beta_1)^3\alpha_2(\alpha_2+\beta_2)^3}$};
\draw (7,-8.5)node{$\frac{1}{\beta_1(\alpha_1+\beta_1)^3\alpha_2(\alpha_2+\beta_2)^3}$};

\draw[step=1cm] (-5,-14) grid (-1,-10);
\draw[step=1cm] (0,-14) grid (4,-10);
\draw[step=1cm] (5,-14) grid (9,-10);
\draw (-4.5,-14.5)node{$t_1$};
\draw (-3.5,-14.5)node{$u_1$};
\draw (-2.5,-14.5)node{$v_1$};
\draw (-1.5,-14.5)node{$w_1$};

\draw (.5,-14.5)node{$t_1$};
\draw (1.5,-14.5)node{$v_1$};
\draw (2.5,-14.5)node{$u_1$};
\draw (3.5,-14.5)node{$w_1$};

\draw (5.5,-14.5)node{$t_1$};
\draw (6.5,-14.5)node{$v_1$};
\draw (7.5,-14.5)node{$w_1$};
\draw (8.5,-14.5)node{$u_1$};
\draw (-5.3,-13.5)node{$v_2$};
\draw (-5.3,-12.5)node{$w_2$};
\draw (-5.3,-11.5)node{$t_2$};
\draw (-5.3,-10.5)node{$u_2$};

\draw (-.3,-13.5)node{$v_2$};
\draw (-.3,-12.5)node{$w_2$};
\draw (-.3,-11.5)node{$t_2$};
\draw (-.3,-10.5)node{$u_2$};

\draw (4.7,-13.5)node{$v_2$};
\draw (4.7,-12.5)node{$w_2$};
\draw (4.7,-11.5)node{$t_2$};
\draw (4.7,-10.5)node{$u_2$};
\draw (-2.5,-13.5)node{$1$};
\draw (-1.5,-12.5)node{$0$};
\draw (-4.5,-11.5)node{$1$};
\draw (-3.5,-10.5)node{$0$};

\draw (1.5,-13.5)node{$1$};
\draw (3.5,-12.5)node{$0$};
\draw (.5,-11.5)node{$1$};
\draw (2.5,-10.5)node{$0$};

\draw (6.5,-13.5)node{$1$};
\draw (7.5,-12.5)node{$0$};
\draw (5.5,-11.5)node{$1$};
\draw (8.5,-10.5)node{$0$};
\draw (-3,-15.5)node{$\frac{1}{\alpha_1^2\beta_2^2N(\alpha+\beta)^2}$};
\draw (2,-15.5)node{$\frac{1}{\alpha_1(\alpha_1+\beta_1)^3\beta_2^2(\alpha_2+\beta_2)^2}$};
\draw (7,-15.5)node{$\frac{1}{\alpha_1(\alpha_1+\beta_1)^3\beta_2^2(\alpha_2+\beta_2)^2}$};
\draw(2,-16.5)node{{\bf{Table (2,1)}}};
\end{tikzpicture}


\end{center}

\begin{center}
\begin{tikzpicture}
\draw[step=1cm] (-5,0) grid (-1,4);
\draw[step=1cm] (0,0) grid (4,4);
\draw[step=1cm] (5,0) grid (9,4);
\draw (-4.5,-.5)node{$v_1$};
\draw (-3.5,-.5)node{$t_1$};
\draw (-2.5,-.5)node{$u_1$};
\draw (-1.5,-.5)node{$w_1$};

\draw (.5,-.5)node{$v_1$};
\draw (1.5,-.5)node{$t_1$};
\draw (2.5,-.5)node{$w_1$};
\draw (3.5,-.5)node{$u_1$};

\draw (5.5,-.5)node{$v_1$};
\draw (6.5,-.5)node{$w_1$};
\draw (7.5,-.5)node{$t_1$};
\draw (8.5,-.5)node{$u_1$};
\draw (-5.3,.5)node{$v_2$};
\draw (-5.3,1.5)node{$t_2$};
\draw (-5.3,2.5)node{$u_2$};
\draw (-5.3,3.5)node{$w_2$};

\draw (-.3,.5)node{$v_2$};
\draw (-.3,1.5)node{$t_2$};
\draw (-.3,2.5)node{$u_2$};
\draw (-.3,3.5)node{$w_2$};

\draw (4.7,.5)node{$v_2$};
\draw (4.7,1.5)node{$t_2$};
\draw (4.7,2.5)node{$u_2$};
\draw (4.7,3.5)node{$w_2$};
\draw (-4.5,.5)node{$1$};
\draw (-2.5,2.5)node{$0$};
\draw (-3.5,1.5)node{$1$};
\draw (-1.5,3.5)node{$0$};

\draw (.5,.5)node{$1$};
\draw (2.5,3.5)node{$0$};
\draw (1.5,1.5)node{$1$};
\draw (3.5,2.5)node{$0$};

\draw (7.5,1.5)node{$1$};
\draw (6.5,3.5)node{$0$};
\draw (5.5,.5)node{$1$};
\draw (8.5,2.5)node{$0$};
\draw (-3,-1.5)node{$\frac{1}{N(\beta)N(\alpha+\beta)^3}$};
\draw (2,-1.5)node{$\frac{1}{N(\beta)N(\alpha+\beta)^3}$};
\draw (7,-1.5)node{$\frac{1}{\alpha_1^2(\alpha_1+\beta_1)^2\alpha_2(\alpha_2+\beta_2)^3}$};
\draw[step=1cm] (-5,-7) grid (-1,-3);
\draw[step=1cm] (0,-7) grid (4,-3);
\draw[step=1cm] (5,-7) grid (9,-3);
\draw (-4.5,-7.5)node{$v_1$};
\draw (-3.5,-7.5)node{$t_1$};
\draw (-2.5,-7.5)node{$u_1$};
\draw (-1.5,-7.5)node{$w_1$};

\draw (.5,-7.5)node{$v_1$};
\draw (1.5,-7.5)node{$t_1$};
\draw (2.5,-7.5)node{$w_1$};
\draw (3.5,-7.5)node{$u_1$};

\draw (5.5,-7.5)node{$v_1$};
\draw (6.5,-7.5)node{$w_1$};
\draw (7.5,-7.5)node{$t_1$};
\draw (8.5,-7.5)node{$u_1$};
\draw (-5.3,-6.5)node{$v_2$};
\draw (-5.3,-5.5)node{$t_2$};
\draw (-5.3,-4.5)node{$w_2$};
\draw (-5.3,-3.5)node{$u_2$};

\draw (-.3,-6.5)node{$v_2$};
\draw (-.3,-5.5)node{$t_2$};
\draw (-.3,-4.5)node{$w_2$};
\draw (-.3,-3.5)node{$u_2$};

\draw (4.7,-6.5)node{$v_2$};
\draw (4.7,-5.5)node{$t_2$};
\draw (4.7,-4.5)node{$w_2$};
\draw (4.7,-3.5)node{$u_2$};
\draw (-4.5,-6.5)node{$1$};
\draw (-3.5,-5.5)node{$1$};
\draw (-1.5,-4.5)node{$0$};
\draw (-2.5,-3.5)node{$0$};

\draw (.5,-6.5)node{$1$};
\draw (1.5,-5.5)node{$1$};
\draw (2.5,-4.5)node{$0$};
\draw (3.5,-3.5)node{$0$};

\draw (5.5,-6.5)node{$1$};
\draw (7.5,-5.5)node{$1$};
\draw (8.5,-3.5)node{$0$};
\draw (6.5,-4.5)node{$0$};
\draw (-3,-8.5)node{$\frac{1}{N(\beta)N(\alpha+\beta)^3}$};
\draw (2,-8.5)node{$\frac{1}{N(\beta)N(\alpha+\beta)^3}$};
\draw (7,-8.5)node{$\frac{1}{\alpha_1^2(\alpha_1+\beta_1)^2\alpha_2^1(\alpha_2+\beta_2)^3}$};

\draw[step=1cm] (-5,-14) grid (-1,-10);
\draw[step=1cm] (0,-14) grid (4,-10);
\draw[step=1cm] (5,-14) grid (9,-10);
\draw (-4.5,-14.5)node{$v_1$};
\draw (-3.5,-14.5)node{$t_1$};
\draw (-2.5,-14.5)node{$u_1$};
\draw (-1.5,-14.5)node{$w_1$};

\draw (.5,-14.5)node{$v_1$};
\draw (1.5,-14.5)node{$t_1$};
\draw (2.5,-14.5)node{$w_1$};
\draw (3.5,-14.5)node{$u_1$};

\draw (5.5,-14.5)node{$v_1$};
\draw (6.5,-14.5)node{$w_1$};
\draw (7.5,-14.5)node{$t_1$};
\draw (8.5,-14.5)node{$u_1$};
\draw (-5.3,-13.5)node{$v_2$};
\draw (-5.3,-12.5)node{$w_2$};
\draw (-5.3,-11.5)node{$t_2$};
\draw (-5.3,-10.5)node{$u_2$};

\draw (-.3,-13.5)node{$v_2$};
\draw (-.3,-12.5)node{$w_2$};
\draw (-.3,-11.5)node{$t_2$};
\draw (-.3,-10.5)node{$u_2$};

\draw (4.7,-13.5)node{$v_2$};
\draw (4.7,-12.5)node{$w_2$};
\draw (4.7,-11.5)node{$t_2$};
\draw (4.7,-10.5)node{$u_2$};
\draw (-4.5,-13.5)node{$1$};
\draw (-3.5,-11.5)node{$1$};
\draw (-1.5,-12.5)node{$0$};
\draw (-2.5,-10.5)node{$0$};

\draw (.5,-13.5)node{$1$};
\draw (1.5,-11.5)node{$1$};
\draw (2.5,-12.5)node{$0$};
\draw (3.5,-10.5)node{$0$};

\draw (5.5,-13.5)node{$1$};
\draw (6.5,-12.5)node{$0$};
\draw (7.5,-11.5)node{$1$};
\draw (8.5,-10.5)node{$0$};
\draw (-3,-15.5)node{$\frac{1}{\alpha_1(\alpha_1+\beta_1)^3\alpha_2^2(\alpha_2+\beta_2)^2}$};
\draw (2,-15.5)node{$\frac{1}{\alpha_1(\alpha_1+\beta_1)^3\alpha_2^2(\alpha_2+\beta_2)^2}$};
\draw (7,-15.5)node{$\frac{1}{N(\beta)^2N(\alpha+\beta)^2}$};
\draw(2,-16.5)node{{\bf{Table (2,2)}}};
\end{tikzpicture}
\end{center}

If $K$ is an imaginary quadratic field we take 
\[C=\N\cup\{\alpha\,|\, Im(\alpha_1)>0\}\]
Then instead of integrals of $f_0$, which is a sum of exponents, we consider an infinite sum of integrals of exponents. In other words, first we consider exponents \[\exp(-\alpha_1t_1-\alpha_2t_2)\times\exp(-\beta_1u_1-\beta_2u_2).\] Then we use iterated integrals over membranes. And finally we sum over element of the cone $C$. Since the sums of such integrals converge for imaginary quadratic fields, we obtain the same shuffle relation as above.


\section{Relations Among Multiple Dedekind Zeta Values}
In this Section, we use integral shuffle relation for the product of Dedekind zeta values at $s=2$ for quadratic fields from Theorem \ref{theorem int shuffle 2,2}.  
From Section 2, we define infinite sum shuffle relation both for real \eqref{eq stuffle 2,2} and for imaginary quadratic fields \eqref{eq stuffle 2,2 C}. 
Using both types of shuffles, we find relations among multiple Dedekind zeta values.

Using the two types of shuffles from Sections 2 and 3, we obtain relations among multiple Dedekind zeta values.

Let $K$ be a real quadratic field. Let ${\cal{O}}_K$ be the ring of integers in $K$. 
If $\alpha$ is in ${\cal{O}}_K$, put $\alpha_1$ and $\alpha_2$ to be its Galois conjugates.
Let 
\[C=\N\{\mu,\nu\}=\{\alpha\in {\cal{O}}_K\,|\, \alpha=a\mu+b\nu,\,a,b\in \N\}.\]
Comparing \eqref{eq int shuffle 2,2} and \eqref{eq stuffle 2,2}, we obtain:
\begin{theorem} Multiple Dedekind zeta values associated to a quadratic number field $K$ satisfy the following relation:
\begin{align}
\nonumber
&2\,_\rho\zeta_{K,C}(2,2;2,2)+\\
\label{eq real 2,2}
&+2\,_{0,1}\zeta_{K,C}(2,2;2,2)+2\,_{1,0}\zeta_{K,C}(2,2;2,2)=\\
\nonumber
=&-\zeta^1_{K;C}(4;4)+8\zeta^1_{K,C}(1,3;1,3)+\\
\nonumber
&+4\zeta^{1}_{K,C}(1,3;2,2)+4\zeta^1_{K,C}(2,2;1,3)+\\
\nonumber
&+2\zeta^{\rho}_{K,C}(2,2;2,2)+8\zeta^\rho_{K,C}(1,3;1,3)+\\
\nonumber
&+4\zeta^{\rho}_{K,C}(1,3;2,2)+4\zeta^\rho_{K,C}(2,2;1,3).
\end{align}
\end{theorem}


Now, let $K$ be an imaginary quadratic field. Let ${\cal{O}}_K$ be the ring of integers in $K$. 
If $\alpha\in K$ let $\alpha_1=\alpha$ and $\alpha_2$ be the complex conjugate of $\alpha$.
Following Kaneko, Gangl, Zagier, \cite{GKZ}, we define
\[C=\N\cup\{\alpha\in {\cal{O}}_K\,|\,Im(\alpha_1)>0\}.\]
From the two types of shuffles for $\zeta_{K,C}(2)\zeta_{K,C}(2)$ from Equations \eqref{eq stuffle 2,2 C} and \eqref{eq int shuffle 2,2}, we obtain: 
\begin{theorem} Multiple Dedekind zeta values associated to an imaginary quadratic number field 
$K$satisfy the following relation:
\begin{align}
\label{eq complex 2,2}
\zeta^1_{K,C}(4,4)-8\zeta^1_{K,C}(1,3;1,3)
=&4\zeta^1(1,3;2,2)+4\zeta^1(2,2;1,3)+\\
\nonumber
&+2\zeta^\rho(2,2;2,2)+2\zeta^\rho(1,3;1,3)+\\
\nonumber
&+4\zeta^\rho(1,3;2,2)+4\zeta^\rho(2,2;1,3).
\end{align}
\end{theorem}



\renewcommand{\em}{\textrm}
\begin{small}

\renewcommand{\refname}{ {\flushleft\normalsize\bf{References}} }
    
\end{small}

\end{document}